\documentclass[12pt]{article}
\usepackage{amssymb, amsmath, amsthm, amscd}
\usepackage[pdftex]{graphics}
\usepackage[utf8]{inputenc}
\usepackage[all,cmtip]{xy}
\usepackage{enumitem}
\usepackage{setspace}
\usepackage{mathdots}

\usepackage[mathscr]{euscript}
\usepackage[colorlinks=true]{hyperref}
\hypersetup{colorlinks   = true}
\hypersetup{linkcolor=blue}
\usepackage{tikz}

\begin{document}

\newtheorem{The}{Theorem}[section]
\newtheorem{Lem}[The]{Lemma}
\newtheorem{Prop}[The]{Proposition}
\newtheorem{Cor}[The]{Corollary}
\newtheorem{Rem}[The]{Remark}
\newtheorem{Obs}[The]{Observation}
\newtheorem{SConj}[The]{Standard Conjecture}
\newtheorem{Titre}[The]{\!\!\!\! }
\newtheorem{Conj}[The]{Conjecture}
\newtheorem{Question}[The]{Question}
\newtheorem{Prob}[The]{Problem}
\newtheorem{Def}[The]{Definition}
\newtheorem{Not}[The]{Notation}
\newtheorem{Claim}[The]{Claim}
\newtheorem{Conc}[The]{Conclusion}
\newtheorem{Ex}[The]{Example}
\newtheorem{Fact}[The]{Fact}
\newtheorem{Formula}[The]{Formula}
\newtheorem{Formulae}[The]{Formulae}
\newtheorem{The-Def}[The]{Theorem and Definition}
\newtheorem{Prop-Def}[The]{Proposition and Definition}
\newtheorem{Lem-Def}[The]{Lemma and Definition}
\newtheorem{Cor-Def}[The]{Corollary and Definition}
\newtheorem{Conc-Def}[The]{Conclusion and Definition}
\newtheorem{Terminology}[The]{Note on terminology}
\newcommand{\C}{\mathbb{C}}
\newcommand{\R}{\mathbb{R}}
\newcommand{\N}{\mathbb{N}}
\newcommand{\Z}{\mathbb{Z}}
\newcommand{\Q}{\mathbb{Q}}
\newcommand{\Proj}{\mathbb{P}}
\newcommand{\Rc}{\mathcal{R}}
\newcommand{\Oc}{\mathcal{O}}
\newcommand{\Vc}{\mathcal{V}}
\newcommand{\Id}{\operatorname{Id}}
\newcommand{\pr}{\operatorname{pr}}
\newcommand{\rk}{\operatorname{rk}}
\newcommand{\del}{\partial}
\newcommand{\delbar}{\bar{\partial}}
\newcommand{\Cdot}{{\raisebox{-0.7ex}[0pt][0pt]{\scalebox{2.0}{$\cdot$}}}}
\newcommand\nilm{\Gamma\backslash G}
\newcommand\frg{{\mathfrak g}}
\newcommand{\fg}{\mathfrak g}
\newcommand{\Oh}{\mathcal{O}}
\newcommand{\Kur}{\operatorname{Kur}}
\newcommand\gc{\frg_\mathbb{C}}
\newcommand\slawek[1]{{\textcolor{red}{#1}}}
\newcommand\dan[1]{{\textcolor{blue}{#1}}}

\begin{center}

{\Large\bf A Moment Map for the Space of Maps to a Balanced Manifold}

\end{center}

\begin{center}

{\large Dan Popovici and Luis Ugarte}

\end{center}

\vspace{1ex}

\noindent{\small{\bf Abstract.} Given a complex balanced manifold $X$ and a compact complex manifold $S$ equipped with a positive volume form $dV>0$ and satisfying an extra condition such that $\mbox{dim}\,S\geq\mbox{dim}\,X -1$, we construct a moment map for the action of the Lie group of biholomorphisms of $S$ that preserve $dV$ onto the space of holomorphic maps $f:S\longrightarrow X$ that satisfy a certain condition with respect to the Bott-Chern cohomology class of the balanced metric of $X$. The purpose is twofold: to study such maps as a possible addition to some very recent hyperbolicity notions involving holomorphic maps with a certain type of growth from some $\C^p$, rather than $S$, to $X$; and to lay the groundwork for a possible future construction of balanced quotients as an analogue of the classical symplectic quotients.}

\vspace{2ex}

\section{Introduction}\label{section:Introduction} Let $X$ be an $n$-dimensional complex manifold. The study of non-constant holomorphic maps $f:\C\longrightarrow X$, if any, has been the subject of a huge literature around the notion of (Brody) hyperbolicity since the foundational work of several authors summed up in Kobayashi's book [Kob70]. Very recently, generalised notions of hyperbolicity were introduced by means of non-degenerate holomorphic maps $f:\C^p\longrightarrow X$, with $p\leq n-1$, satisfying certain growth conditions (cf. e.g. [MP22a], [MP22b], [KP23]).

In this paper, we take a related but slightly different approach to holomorphic maps into $X$ by supposing that their origin is a compact complex manifold $S$ of dimension $d\geq n-1$. These maps are assumed to satisfy a compatibility condition with a {\it balanced} structure that is supposed to exist on $X$. Thus, the set-up is, to some extent, a complex balanced analogue of the classical symplectic situation described, for example, in [Don99].


To our knowledge, a few attempts have already been made in the literature to develop a moment map theory in the Hermitian setting of possibly non-K\"ahler special metrics. 
For instance, in [GRT23], for a compact complex manifold $X$ with a smooth volume form compatible with the orientation and for any $l\in\R\setminus\{2\}$, a moment map $\mu_l$ is defined on $T\Omega^{1,1}_{>0}$ so that the zeros of $\mu_l$ are given by the $l$-conformally balanced metrics on $X$. Here $T\Omega^{1,1}_{>0}=\Omega^{1,1}_{>0} \times \Omega^{1,1}_{\R}$, where $\Omega^{1,1}_{>0}\subset \Omega^{1,1}_{\R}$ is the space of real positive (1,1)-forms on $X$ (see [GRT23, Proposition 5.4]). This moment map interpretation of the conformally balanced metrics is further developed in [GRT23, Proposition 5.14] to characterise the solutions of the Calabi system, a coupled system of equations unifying the classical Calabi problem and the Hull-Strominger system, as a moment map condition.

In the framework of locally conformally K\"ahler geometry, given a 
compact manifold $M$ endowed with a locally conformally
symplectic structure $\omega$, a moment map $\mu$ is defined in [ACPS23] on the space ${\mathcal J}(\omega)^K$ that consists of $K$-invariant, compatible integrable almost complex structures on $(M, \omega)$. Here $K$ is a compact, connected Lie group  acting effectively on $M$ and satisfying certain additional conditions in terms of $\omega$. This is then used to extend the Fujiki-Donaldson interpretation of the scalar curvature as a moment map in K\"ahler geometry to locally conformally K\"ahler manifolds (see [ACPS23, Theorem 4.3]).

The approach we take in this paper consists in considering a moment map defined on a certain space of maps rather than on a space of forms, in the spirit of [Don99]. We hope that this approach might be further developed to construct ``balanced quotients" and possibly new examples of balanced manifolds in future work.

\vspace{2ex}

\noindent {\bf Acknowledgments.} The authors are very grateful to C. Scarpa for useful comments and for sharing his thoughts about the context of moment maps and their applicability to various problems in complex geometry.

\section{Definition and well-definedness of the moment map}\label{section:moment-map_w-d} Let $X$ be an $n$-dimensional (not necessarily compact) complex manifold supposed to carry {\bf balanced} metrics. We fix such a metric $\omega$. This is a Hermitian metric (identified with a $C^\infty$ positive definite $(1,\,1)$-form $\omega$ on $X$) satisfying the extra condition $d\omega^{n-1}=0$. In particular, $\omega_{n-1}:=\omega^{n-1}/(n-1)!$ defines a Bott-Chern cohomology class $[\omega_{n-1}]_{BC}\in H^{n-1,\,n-1}_{BC}(X,\,\C)$ on $X$.  

On the other hand, let $S$ be a compact complex manifold of dimension $d\geq n-1$ on which a volume form $dV>0$ has been fixed. This means that $dV$ is a $C^\infty$ $(d,\,d)$-form on $S$ that is positive at every point. We consider the following space of maps: \begin{eqnarray*}{\mathscr X}:=\bigg\{f:S\longrightarrow X\,\mid\,f\hspace{1ex}\mbox{is holomorphic and}\hspace{1ex} f^\star[\omega_{n-1}]_{BC} = 0\bigg\},\end{eqnarray*} where $f^\star[\omega_{n-1}]_{BC} = [f^\star\omega_{n-1}]_{BC}\in  H^{n-1,\,n-1}_{BC}(S,\,\C)$ is the Bott-Chern cohomology class of the pullback (a smooth $(n-1,\,n-1)$-form) under $f$ of $\omega_{n-1}$ to $S$.

\vspace{2ex} 

$\bullet$ Note that if $d=n-1$ and there exists a point $x\in S$ such that $f$ is non-degenerate at $x$ (in the sense that the differential map $d_xf:T^{1,\,0}_xS\longrightarrow T^{1,\,0}_{f(x)}X$ is of maximal rank), the condition $f^\star[\omega_{n-1}]_{BC} = 0$ is impossible. Indeed, in that case, $f^\star\omega_{n-1}$ would be a smooth form of top bidegree $(n-1,\,n-1)$ on $S$ lying in the image of $\partial\bar\partial$. By Stokes and the compactness of $S$, this would lead to $\int_Sf^\star\omega_{n-1} = 0$, which would imply, thanks to $f^\star\omega_{n-1}\geq 0$ everywhere on $S$, that $f^\star\omega_{n-1} = 0$, contradicting the strict positivity of $f^\star\omega_{n-1}$ at $x$.

\vspace{1ex}

$\bullet$ The case $d=n$ is especially interesting to us since the condition $f^\star[\omega_{n-1}]_{BC} = 0$ is then equivalent to the smooth semi-positive definite $(n-1,\,n-1)$-form $f^\star\omega_{n-1}$ being $\partial\bar\partial$-exact on $S$. When $f$ is non-degenerate at every point of $S$, this implies that $f^\star\omega_{n-1}$ is a {\it degenerate balanced metric} on $S$, namely a smooth, positive definite, $d$-exact $(n-1,\,n-1)$-form on the $n$-dimensional complex manifold $S$. (See [Pop15, Proposition 5.4.] for equivalent descriptions of compact complex manifolds admitting degenerate balanced metrics and [MP22a, $\S2.3$] for a discussion of these manifolds as special cases of {\it balanced hyperbolic manifolds} introduced therein and a reminder of the two classes of such manifolds known so far, originating in [Fri89] and [Yac98].)

\vspace{1ex}

$\bullet$ The case $d\geq n+1$ seems equally worthy of further investigation since then the condition $f^\star[\omega_{n-1}]_{BC} = 0$ yields a $\partial\bar\partial$-exact, smooth $(n-1,\,n-1)$-form $f^\star\omega_{n-1}$ on $S$ generalising the notion of $(n-1)$-K\"ahler structure on $S$ introduced in [AA87].

\vspace{2ex}

We can think of ${\mathscr X}$ as a kind of ``complex manifold'' (although it may not have and we will not suppose it to have such a structure) by defining its ``tangent space'' of type $(1,\,0)$ at every point $f\in{\mathscr X}$ as the space of global holomorphic sections on $S$ of the pullback bundle $f^\star T^{1,\,0}X\subset T^{1,\,0}S$: \begin{eqnarray*}T^{1,\,0}_f{\mathscr X}:=H^0(S,\,f^\star T^{1,\,0}X).\end{eqnarray*} In particular, the tangent vectors to ${\mathscr X}$ are holomorphic vector fields of type $(1,\,0)$ on $S$.

The space ${\mathscr X}$ carries a natural $C^\infty$ weakly strictly positive $(n-1,\,n-1)$-form $\Omega$ defined, for every $f\in{\mathscr X}$, by \begin{eqnarray}\label{eqn:Omega_def}\Omega(f)(v_1,\dots , v_{n-1},\,\bar{w}_1,\dots , \bar{w}_{n-1}):=\int\limits_S(f^\star\omega_{n-1})(v_1,\dots , v_{n-1},\,\bar{w}_1,\dots , \bar{w}_{n-1})\,dV\end{eqnarray} for any sections $v_1,\dots , v_{n-1},\,w_1,\dots , w_{n-1}\in T^{1,\,0}_f{\mathscr X}:=H^0(S,\,f^\star T^{1,\,0}X)$. Note that $d\Omega=0$ on ${\mathscr X}$ since $d(f^\star\omega_{n-1})=0$ on $S$, so $\Omega$ is an $(n-1)$-K\"ahler structure on ${\mathscr X}$ (assuming that ${\mathscr X}$ is a complex manifold -- otherwise, the $(n-1)$-K\"ahler structure can be construed formally).

  We also consider the Lie group of biholomorphisms of $S$ that preserve the volume form $dV$: \begin{eqnarray*}{\cal G}:= \bigg\{\varphi:S\longrightarrow S \mid\,\varphi\hspace{1ex}\mbox{is a biholomorphism and}\hspace{1ex}\varphi^\star(dV) = dV\bigg\}.\end{eqnarray*} The Lie group ${\cal G}$ acts on ${\mathscr X}$ by composition on the right and this action preserves $\Omega$. It is thus legitimate to seek a moment map for this action.

 The Lie algebra of ${\cal G}$, viewed as its tangent space at the identity map, consists of holomorphic vector fields $\xi$ of type $(1,\,0)$ on $S$ that preserve the volume form $dV$ in the sense that the $(1,\,0)$-Lie derivative (see Definition \ref{Def:Lie-deriv_def}) w.r.t. $\xi$ of $dV$ vanishes: \begin{eqnarray*}\frg = Lie({\cal G}) = T^{1,\,0}_{\mbox{Id}}{\cal G} = \bigg\{\xi\in H^0(S,\,T^{1,\,0}S)\,\mid\,L^{1,\,0}_\xi(dV) = 0\bigg\}.\end{eqnarray*}

We define the following vector subspace of $Lie({\cal G})^{n-2}\times\overline{Lie({\cal G})}^{n-2}$: \begin{eqnarray*}{\cal P}_\frg:=\bigg\{(\xi_1,\dots , \xi_{n-2},\,\bar\eta_1,\dots , \bar\eta_{n-2})\in Lie({\cal G})^{n-2}\times\overline{Lie({\cal G})}^{n-2}\,\mid\,(\ref{eqn:P_g_def_1})\hspace{1ex}\mbox{and}\hspace{1ex}(\ref{eqn:P_g_def_2})\hspace{1ex}\mbox{hold}\bigg\},\end{eqnarray*} where the conditions (\ref{eqn:P_g_def_1}) and (\ref{eqn:P_g_def_2}) imposed on $(\xi_1,\dots , \xi_{n-2},\,\bar\eta_1,\dots , \bar\eta_{n-2})$ are the following: \begin{eqnarray}\label{eqn:def_cal-P_1}\nonumber\bigg[\sum\limits_{l=2}^{n-1}\sum\limits_{r=1}^{n-l-1}(-1)^{n+1-l-r}[\bar\eta_{n-l},\,\bar\eta_r]\lrcorner\bar\eta_{n-2}\lrcorner\dots\lrcorner\widehat{\bar\eta_{n-l}}\lrcorner\dots\lrcorner\widehat{\bar\eta_r}\lrcorner\dots\lrcorner\bar\eta_1\lrcorner\xi_{n-2}\lrcorner\dots\lrcorner\xi_1\end{eqnarray}\begin{eqnarray}\label{eqn:P_g_def_1} + \sum\limits_{l=2}^{n-1}\sum\limits_{r=1}^{n-2}(-1)^{l+r+1}[\bar\eta_{n-l},\,\xi_r]\lrcorner\bar\eta_{n-2}\lrcorner\dots\lrcorner\widehat{\bar\eta_{n-l}}\lrcorner\dots\lrcorner\bar\eta_1\lrcorner\xi_{n-2}\lrcorner\dots\lrcorner\widehat{\xi_r}\lrcorner\dots\lrcorner\xi_1\bigg]\lrcorner dV = 0\end{eqnarray}

    \noindent and \begin{eqnarray}\label{eqn:def_cal-P_2}\nonumber\bigg[\sum\limits_{l=2}^{n-1}\sum\limits_{r=1}^{n-l-1}(-1)^{n+1-l-r}[\xi_{n-l},\,\xi_r]\lrcorner\xi_{n-2}\lrcorner\dots\lrcorner\widehat{\xi_{n-l}}\lrcorner\dots\lrcorner\widehat{\xi_r}\lrcorner\dots\lrcorner\xi_1\lrcorner\bar\eta_{n-2}\lrcorner\dots\lrcorner\bar\eta_1\end{eqnarray}\begin{eqnarray}\label{eqn:P_g_def_2} + \sum\limits_{l=2}^{n-1}\sum\limits_{r=1}^{n-2}(-1)^{l+r+1}[\xi_{n-l},\,\bar\eta_r]\lrcorner\xi_{n-2}\lrcorner\dots\lrcorner\widehat{\xi_{n-l}}\lrcorner\dots\lrcorner\xi_1\lrcorner\bar\eta_{n-2}\lrcorner\dots\lrcorner\widehat{\bar\eta_r}\lrcorner\dots\lrcorner\bar\eta_1\bigg]\lrcorner dV = 0.\end{eqnarray}

      \vspace{2ex}

      Now, for every $f\in{\mathscr X}$, $f^\star\omega_{n-1}$ is a smooth real $\partial\bar\partial$-exact $(n-1,\,n-1)$-form on $S$, hence there exists a real form $\Gamma_f\in C^\infty_{n-2,\,n-2}(S,\,\C)$ such that \begin{eqnarray}\label{eqn:f-omega_n1_potential}f^\star\omega_{n-1} = i\partial\bar\partial\Gamma_f \hspace{3ex} \mbox{on}\hspace{1ex} S.\end{eqnarray}

\begin{Def}\label{Def:moment-map_def} With every $n$-dimensional complex {\bf balanced} manifold $(X,\,\omega_{n-1})$ and every $d$-dimensional {\bf compact} complex manifold equipped with a positive {\bf volume form} $(S,\,dV)$ such that $d\geq n-1$, we associate the map $\mu:{\mathscr X}\longrightarrow{\cal P}_\frg^\star$ defined by \begin{eqnarray}\label{eqn:moment-map_def}\mu(f)(\xi_1,\dots , \xi_{n-2},\,\bar\eta_1,\dots , \bar\eta_{n-2}):=i\int\limits_S\Gamma_f(\xi_1,\dots , \xi_{n-2},\,\bar\eta_1,\dots , \bar\eta_{n-2})\,dV\end{eqnarray} for every $f\in{\mathscr X}$ and every $(\xi_1,\dots , \xi_{n-2},\,\bar\eta_1,\dots , \bar\eta_{n-2})\in{\cal P}_\frg$.

\end{Def}

Our first task is to prove the following

\begin{Prop}\label{Prop:moment-map_well-defined} Suppose that $H^{n-2,\,n-2}_A(S,\,\C)=\{0\}$. Then, the map $\mu$ of Definition \ref{Def:moment-map_def} is {\bf well defined} in the sense that the expression on the right-hand side of (\ref{eqn:moment-map_def}) is independent of the choice of $\Gamma_f$ satisfying property (\ref{eqn:f-omega_n1_potential}). 

\end{Prop}

\noindent {\it Proof.} Suppose $\Gamma_1$ and $\Gamma_2$ are real $(n-2,\,n-2)$-forms on $S$ such that $f^\star\omega_{n-1} = i\partial\bar\partial\Gamma_1 = i\partial\bar\partial\Gamma_2$. Then, $\Gamma_1 - \Gamma_2\in\ker(\partial\bar\partial)$ and we get an Aeppli cohomology class $[\Gamma_1 - \Gamma_2]_A\in H^{n-2,\,n-2}_A(S,\,\C)$. This cohomology space being zero, by hypothesis, and the form $\Gamma_1 - \Gamma_2$ being real, there exists a smooth $(n-2,\,n-3)$-form $\beta$ such that \begin{eqnarray*}\Gamma_1 - \Gamma_2 = \partial\bar\beta + \bar\partial\beta \hspace{3ex} \mbox{on}\hspace{1ex} S.\end{eqnarray*}

To prove the contention, we need to show that \begin{eqnarray}\label{eqn:moment-map_w-d_proof_1}\int\limits_S(\partial\bar\beta + \bar\partial\beta)(\xi_1,\dots , \xi_{n-2},\,\bar\eta_1,\dots , \bar\eta_{n-2})\,dV = 0\end{eqnarray} for every $(\xi_1,\dots , \xi_{n-2},\,\bar\eta_1,\dots , \bar\eta_{n-2})\in{\cal P}_\frg$.

  \begin{Claim}\label{Claim:del-beta_bar_parts-0} For every $(\xi_1,\dots , \xi_{n-2},\,\bar\eta_1,\dots , \bar\eta_{n-2})\in{\cal P}_\frg$, the following equalities hold: \begin{eqnarray}\label{eqn:del-beta_bar_parts}\nonumber(\partial\bar\beta)(\xi_1,\dots , \xi_{n-2},\,\bar\eta_1,\dots , \bar\eta_{n-2})\,dV & = & (-1)^{n-2}(\partial\bar\beta)\wedge(\xi_1\lrcorner\dots\lrcorner\xi_{n-2}\lrcorner\bar\eta_1\lrcorner\dots\lrcorner\bar\eta_{n-2}\lrcorner dV) \\
 (\bar\partial\beta)(\xi_1,\dots , \xi_{n-2},\,\bar\eta_1,\dots , \bar\eta_{n-2})\,dV & = & (-1)^{n-2}(\bar\partial\beta)\wedge(\xi_1\lrcorner\dots\lrcorner\xi_{n-2}\lrcorner\bar\eta_1\lrcorner\dots\lrcorner\bar\eta_{n-2}\lrcorner dV).\end{eqnarray}

\end{Claim}

  \noindent {\it Proof of Claim.} It suffices to prove the first equality in (\ref{eqn:del-beta_bar_parts}). Applying repeatedly the definition of the contraction of a differential form by a vector field, we get:  \begin{eqnarray*}(\partial\bar\beta)(\xi_1,\dots , \xi_{n-2},\,\bar\eta_1,\dots , \bar\eta_{n-2})\,dV = \bigg(\bar\eta_{n-2}\lrcorner\dots\lrcorner\bar\eta_1\lrcorner\xi_{n-2}\lrcorner\dots\lrcorner\xi_1\lrcorner(\partial\bar\beta)\bigg)\,dV:=T_1.\end{eqnarray*} (So, we denote by $T_1$ the second quantity above.) We will compute $T_1$ in successive stages.

\vspace{1ex}  

$\bullet$ Starting from the equality $(\bar\eta_{n-3}\lrcorner\dots\lrcorner\bar\eta_1\lrcorner\xi_{n-2}\lrcorner\dots\lrcorner\xi_1\lrcorner(\partial\bar\beta))\wedge dV = 0$, which holds trivially since the form on the left has bidegree $(d,\,d+1)$ on the $d$-dimensional complex manifold $S$, we trivially infer that \begin{eqnarray*}\bar\eta_{n-2}\lrcorner\bigg[\bigg(\bar\eta_{n-3}\lrcorner\dots\lrcorner\bar\eta_1\lrcorner\xi_{n-2}\lrcorner\dots\lrcorner\xi_1\lrcorner(\partial\bar\beta)\bigg)\wedge dV\bigg]  = 0.\end{eqnarray*} From this, we get: \begin{eqnarray*}T_1 = \bigg(\bar\eta_{n-3}\lrcorner\dots\lrcorner\bar\eta_1\lrcorner\xi_{n-2}\lrcorner\dots\lrcorner\xi_1\lrcorner(\partial\bar\beta)\bigg)\wedge(\bar\eta_{n-2}\lrcorner dV):=S_1,\end{eqnarray*} where $S_1$ is the name we give to the second quantity.

\vspace{1ex}  

$\bullet$ Starting from the equality $(\bar\eta_{n-4}\lrcorner\dots\lrcorner\bar\eta_1\lrcorner\xi_{n-2}\lrcorner\dots\lrcorner\xi_1\lrcorner(\partial\bar\beta))\wedge(\bar\eta_{n-2}\lrcorner dV) = 0$, which holds trivially since the form on the left has bidegree $(d,\,d+1)$ on the $d$-dimensional complex manifold $S$, we trivially infer that \begin{eqnarray*}\bar\eta_{n-3}\lrcorner\bigg[\bigg(\bar\eta_{n-4}\lrcorner\dots\lrcorner\bar\eta_1\lrcorner\xi_{n-2}\lrcorner\dots\lrcorner\xi_1\lrcorner(\partial\bar\beta)\bigg)\wedge(\bar\eta_{n-2}\lrcorner dV)\bigg]  = 0.\end{eqnarray*} From this, we get: \begin{eqnarray*}S_1 = - \bigg(\bar\eta_{n-4}\lrcorner\dots\lrcorner\bar\eta_1\lrcorner\xi_{n-2}\lrcorner\dots\lrcorner\xi_1\lrcorner(\partial\bar\beta)\bigg)\wedge(\bar\eta_{n-3}\lrcorner\bar\eta_{n-2}\lrcorner dV) := -S_2,\end{eqnarray*} where $-S_2$ is the name we give to the second quantity.

\vspace{1ex}  

$\bullet$ Starting from the equality $(\bar\eta_{n-5}\lrcorner\dots\lrcorner\bar\eta_1\lrcorner\xi_{n-2}\lrcorner\dots\lrcorner\xi_1\lrcorner(\partial\bar\beta))\wedge(\bar\eta_{n-3}\lrcorner\bar\eta_{n-2}\lrcorner dV) = 0$, which holds trivially since the form on the left has bidegree $(d,\,d+1)$ on the $d$-dimensional complex manifold $S$, we trivially infer that \begin{eqnarray*}\bar\eta_{n-4}\lrcorner\bigg[\bigg(\bar\eta_{n-5}\lrcorner\dots\lrcorner\bar\eta_1\lrcorner\xi_{n-2}\lrcorner\dots\lrcorner\xi_1\lrcorner(\partial\bar\beta)\bigg)\wedge(\bar\eta_{n-3}\lrcorner\bar\eta_{n-2}\lrcorner dV)\bigg]  = 0.\end{eqnarray*} From this, we get: \begin{eqnarray*}S_2 = \bigg(\bar\eta_{n-5}\lrcorner\dots\lrcorner\bar\eta_1\lrcorner\xi_{n-2}\lrcorner\dots\lrcorner\xi_1\lrcorner(\partial\bar\beta)\bigg)\wedge(\bar\eta_{n-4}\lrcorner\bar\eta_{n-3}\lrcorner\bar\eta_{n-2}\lrcorner dV) := S_3,\end{eqnarray*} where $S_3$ is the name we give to the second quantity.

\vspace{1ex}  

$\bullet$ Continuing in this way, by induction, we get: $$T_1 = S_1 = -S_2 = -S_3 = S_4 = S_5 = \dots = (-1)^{n-2}\,S_{2(n-2)},$$ where, for every $k$, $S_k$ is defined in a way similar to $S_1, S_2, S_3$ with $k$ forms, taken in the order specified above from the $\bar\eta_l$'s and the $\xi_j$'s, contracting $dV$ in the second factor of the exterior product. This proves the first equality in (\ref{eqn:del-beta_bar_parts}).  The second one is proved in the same way after replacing $\partial\bar\beta$ by $\bar\partial\beta$. \hfill $\Box$

\vspace{2ex}

\noindent {\it Sequel to the proof of Proposition \ref{Prop:moment-map_well-defined}.} Recall that we have to prove equality (\ref{eqn:moment-map_w-d_proof_1}). Thanks to Claim \ref{Claim:del-beta_bar_parts-0}, for any fixed $(\xi_1,\dots , \xi_{n-2},\,\bar\eta_1,\dots , \bar\eta_{n-2})\in{\cal P}_\frg$, the integral of (\ref{eqn:moment-map_w-d_proof_1}) multiplied by $(-1)^{n-2}$ is given by the first equality below: \begin{eqnarray}\label{eqn:moment-map_w-d_proof_2}\nonumber & & (-1)^{n-2}\,\int\limits_S(\partial\bar\beta + \bar\partial\beta)(\xi_1,\dots , \xi_{n-2},\,\bar\eta_1,\dots , \bar\eta_{n-2})\,dV \\
\nonumber  & = & \int\limits_S(\partial\bar\beta)\wedge(\xi_1\lrcorner\dots\lrcorner\xi_{n-2}\lrcorner\bar\eta_1\lrcorner\dots\lrcorner\bar\eta_{n-2}\lrcorner dV) + \int\limits_S(\bar\partial\beta)\wedge(\xi_1\lrcorner\dots\lrcorner\xi_{n-2}\lrcorner\bar\eta_1\lrcorner\dots\lrcorner\bar\eta_{n-2}\lrcorner dV) \\
\nonumber  & = & \int\limits_S\partial\bigg[\bar\beta\wedge(\xi_1\lrcorner\dots\lrcorner\xi_{n-2}\lrcorner\bar\eta_1\lrcorner\dots\lrcorner\bar\eta_{n-2}\lrcorner dV)\bigg] + \int\limits_S\bar\beta\wedge\partial(\xi_1\lrcorner\dots\lrcorner\xi_{n-2}\lrcorner\bar\eta_1\lrcorner\dots\lrcorner\bar\eta_{n-2}\lrcorner dV) \\
\nonumber  & + & \int\limits_S\bar\partial\bigg[\beta\wedge(\xi_1\lrcorner\dots\lrcorner\xi_{n-2}\lrcorner\bar\eta_1\lrcorner\dots\lrcorner\bar\eta_{n-2}\lrcorner dV)\bigg] + \int\limits_S\beta\wedge\bar\partial(\xi_1\lrcorner\dots\lrcorner\xi_{n-2}\lrcorner\bar\eta_1\lrcorner\dots\lrcorner\bar\eta_{n-2}\lrcorner dV) \\
& = & \int\limits_S\bar\beta\wedge\partial(\xi_1\lrcorner\dots\lrcorner\xi_{n-2}\lrcorner\bar\eta_1\lrcorner\dots\lrcorner\bar\eta_{n-2}\lrcorner dV) + \int\limits_S\beta\wedge\bar\partial(\xi_1\lrcorner\dots\lrcorner\xi_{n-2}\lrcorner\bar\eta_1\lrcorner\dots\lrcorner\bar\eta_{n-2}\lrcorner dV),\end{eqnarray} where the last equality follows from Stokes's theorem thanks to $S$ being compact.

To continue, we need to compute the integrands on the last line of (\ref{eqn:moment-map_w-d_proof_2}). To this end, we will prove the following

\begin{Claim}\label{Claim:del-beta_bar_parts} For every $(\xi_1,\dots , \xi_{n-2},\,\bar\eta_1,\dots , \bar\eta_{n-2})\in Lie({\cal G})^{n-2}\times\overline{Lie({\cal G})}^{n-2}$, the following equalities hold: \begin{eqnarray}\label{eqn:integrands_w-d_1}\nonumber & & \partial(\xi_1\lrcorner\dots\lrcorner\xi_{n-2}\lrcorner\bar\eta_1\lrcorner\dots\lrcorner\bar\eta_{n-2}\lrcorner dV) \\
 \nonumber   & = & \bigg[\sum\limits_{l=2}^{n-1}\sum\limits_{r=1}^{n-l-1}(-1)^{n+1-l-r}[\xi_{n-l},\,\xi_r]\lrcorner\xi_{n-2}\lrcorner\dots\lrcorner\widehat{\xi_{n-l}}\lrcorner\dots\lrcorner\widehat{\xi_r}\lrcorner\dots\lrcorner\xi_1\lrcorner\bar\eta_{n-2}\lrcorner\dots\lrcorner\bar\eta_1 \\
  & + & \sum\limits_{l=2}^{n-1}\sum\limits_{r=1}^{n-2}(-1)^{l+r+1}[\xi_{n-l},\,\bar\eta_r]\lrcorner\xi_{n-2}\lrcorner\dots\lrcorner\widehat{\xi_{n-l}}\lrcorner\dots\lrcorner\xi_1\lrcorner\bar\eta_{n-2}\lrcorner\dots\lrcorner\widehat{\bar\eta_r}\lrcorner\dots\lrcorner\bar\eta_1\bigg]\lrcorner dV  \end{eqnarray}

\noindent and

\begin{eqnarray}\label{eqn:integrands_w-d_2}\nonumber & & \bar\partial(\xi_1\lrcorner\dots\lrcorner\xi_{n-2}\lrcorner\bar\eta_1\lrcorner\dots\lrcorner\bar\eta_{n-2}\lrcorner dV) \\
 \nonumber & = & \bigg[\sum\limits_{l=2}^{n-1}\sum\limits_{r=1}^{n-l-1}(-1)^{n+1-l-r}[\bar\eta_{n-l},\,\bar\eta_r]\lrcorner\bar\eta_{n-2}\lrcorner\dots\lrcorner\widehat{\bar\eta_{n-l}}\lrcorner\dots\lrcorner\widehat{\bar\eta_r}\lrcorner\dots\lrcorner\bar\eta_1\lrcorner\xi_{n-2}\lrcorner\dots\lrcorner\xi_1 \\
 & + & \sum\limits_{l=2}^{n-1}\sum\limits_{r=1}^{n-2}(-1)^{l+r+1}[\bar\eta_{n-l},\,\xi_r]\lrcorner\bar\eta_{n-2}\lrcorner\dots\lrcorner\widehat{\bar\eta_{n-l}}\lrcorner\dots\lrcorner\bar\eta_1\lrcorner\xi_{n-2}\lrcorner\dots\lrcorner\widehat{\xi_r}\lrcorner\dots\lrcorner\xi_1\bigg]\lrcorner dV.\end{eqnarray}

\end{Claim}

\noindent {\it Proof of Claim.} We will prove the second equality in (\ref{eqn:integrands_w-d_2}). The first one can be proved in a similar fashion. Note that \begin{eqnarray*}\bar\partial(\xi_1\lrcorner\dots\lrcorner\xi_{n-2}\lrcorner\bar\eta_1\lrcorner\dots\lrcorner\bar\eta_{n-2}\lrcorner dV) = (-1)^{n-2}\,\bar\partial(
  \bar\eta_{n-2}\lrcorner\dots\lrcorner\bar\eta_1\lrcorner\xi_{n-2}\lrcorner\dots\lrcorner\xi_1\lrcorner dV):=(-1)^{n-2}\,T.\end{eqnarray*}

We will compute $T$ using the properties of the Lie derivatives of types $(1,\,0)$ and $(0,\,1)$ listed and proved in $\S.$\ref{section:Lie-derivatives}. We get: \begin{eqnarray*}T & = & \bar\partial(\bar\eta_{n-2}\lrcorner\dots\lrcorner\bar\eta_1\lrcorner\xi_{n-2}\lrcorner\dots\lrcorner\xi_1\lrcorner dV) \\
    & = & L^{0,\,1}_{\bar\eta_{n-2}}(\bar\eta_{n-3}\lrcorner\dots\lrcorner\bar\eta_1\lrcorner\xi_{n-2}\lrcorner\dots\lrcorner\xi_1\lrcorner dV) - \bar\eta_{n-2}\lrcorner\bar\partial(\bar\eta_{n-3}\lrcorner\dots\lrcorner\bar\eta_1\lrcorner\xi_{n-2}\lrcorner\dots\lrcorner\xi_1\lrcorner dV) \\
   & = & [\bar\eta_{n-2},\,\bar\eta_{n-3}]\lrcorner(\bar\eta_{n-4}\lrcorner\dots\lrcorner\bar\eta_1\lrcorner\xi_{n-2}\lrcorner\dots\lrcorner\xi_1\lrcorner dV) + \bar\eta_{n-3}\lrcorner L^{0,\,1}_{\bar\eta_{n-2}}(\bar\eta_{n-4}\lrcorner\dots\lrcorner\bar\eta_1\lrcorner\xi_{n-2}\lrcorner\dots\lrcorner\xi_1\lrcorner dV) \\
  & - & \bar\eta_{n-2}\lrcorner L^{0,\,1}_{\bar\eta_{n-3}}(\bar\eta_{n-4}\lrcorner\dots\lrcorner\bar\eta_1\lrcorner\xi_{n-2}\lrcorner\dots\lrcorner\xi_1\lrcorner dV) + \bar\eta_{n-2}\lrcorner\bar\eta_{n-3}\lrcorner\bar\partial(\bar\eta_{n-4}\lrcorner\dots\lrcorner\bar\eta_1\lrcorner\xi_{n-2}\lrcorner\dots\lrcorner\xi_1\lrcorner dV),\end{eqnarray*} where definition (\ref{eqn:Lie-deriv_0-1_def}) was used for $L^{0,\,1}_{\bar\eta_{n-2}}$ to get the first equality, while property (\ref{eqn:L_commutations-1}) and again definition (\ref{eqn:Lie-deriv_0-1_def}) were used to get the second equality.

The contention follows by induction on the number of terms by repeating the above arguments, starting from the equalities: \begin{eqnarray*}L^{0,\,1}_{\bar\eta_{n-2}}(\bar\eta_{n-4}\lrcorner\dots\lrcorner\bar\eta_1\lrcorner\xi_{n-2}\lrcorner\dots\lrcorner\xi_1\lrcorner dV) & = & \bar\eta_{n-4}\lrcorner L^{0,\,1}_{\bar\eta_{n-2}}(\bar\eta_{n-5}\lrcorner\dots\lrcorner\bar\eta_1\lrcorner\xi_{n-2}\lrcorner\dots\lrcorner\xi_1\lrcorner dV) \\
  & + & [\bar\eta_{n-2},\,\bar\eta_{n-4}]\lrcorner(\bar\eta_{n-5}\lrcorner\dots\lrcorner\bar\eta_1\lrcorner\xi_{n-2}\lrcorner\dots\lrcorner\xi_1\lrcorner dV),\end{eqnarray*} its analogue expressing $L^{0,\,1}_{\bar\eta_{n-3}}(\bar\eta_{n-4}\lrcorner\dots\lrcorner\bar\eta_1\lrcorner\xi_{n-2}\lrcorner\dots\lrcorner\xi_1\lrcorner dV)$ and \begin{eqnarray*}\bar\partial(\bar\eta_{n-4}\lrcorner\dots\lrcorner\bar\eta_1\lrcorner\xi_{n-2}\lrcorner\dots\lrcorner\xi_1\lrcorner dV) & = & L^{0,\,1}_{\bar\eta_{n-4}}(\bar\eta_{n-5}\lrcorner\dots\lrcorner\bar\eta_1\lrcorner\xi_{n-2}\lrcorner\dots\lrcorner\xi_1\lrcorner dV) \\
  & - & \bar\eta_{n-4}\lrcorner\bar\partial(\bar\eta_{n-5}\lrcorner\dots\lrcorner\bar\eta_1\lrcorner\xi_{n-2}\lrcorner\dots\lrcorner\xi_1\lrcorner dV),\end{eqnarray*} the first two of which following from property (\ref{eqn:L_commutations-1}), while the third follows from definition (\ref{eqn:Lie-deriv_0-1_def}).

Let us only mention that when a $(0,\,1)$-Lie derivative reaches a contraction by some $\xi_j$, we can apply property (vi) of Lemma \ref{Lem:Lie-deriv-prop_0-1} to write, for example, \begin{eqnarray*}\bar\eta_{n-2}\lrcorner\bar\eta_{n-4}\lrcorner\dots\lrcorner\bar\eta_1\lrcorner L^{0,\,1}_{\bar\eta_{n-3}}(\xi_{n-2}\lrcorner\dots\lrcorner\xi_1\lrcorner dV) & = & \bar\eta_{n-2}\lrcorner\bar\eta_{n-4}\lrcorner\dots\lrcorner\bar\eta_1\lrcorner\xi_{n-2}\lrcorner L^{0,\,1}_{\bar\eta_{n-3}}(\xi_{n-3}\lrcorner\dots\lrcorner\xi_1\lrcorner dV) \\
 & + & \bar\eta_{n-2}\lrcorner\bar\eta_{n-4}\lrcorner\dots\lrcorner\bar\eta_1\lrcorner[\bar\eta_{n-3},\,\xi_{n-2}]\lrcorner\xi_{n-3}\lrcorner\dots\lrcorner\xi_1\lrcorner dV)\end{eqnarray*} since $\eta_{n-3}$ is holomorphic. Moreover, $L^{0,\,1}_{\bar\eta_k}(dV)=0$ for each $k$ since $\eta_k\in Lie({\cal G})$, so $L^{1,\,0}_{\eta_k}(dV)=0$.  \hfill $\Box$

\vspace{2ex}

\noindent {\it End of proof of Proposition \ref{Prop:moment-map_well-defined}.} Formulae (\ref{eqn:integrands_w-d_1}) and (\ref{eqn:integrands_w-d_2}), together with the definition of ${\cal P}_\frg$ (see (\ref{eqn:def_cal-P_1}) and (\ref{eqn:def_cal-P_2})), show that the hypothesis $(\xi_1,\dots , \xi_{n-2},\,\bar\eta_1,\dots , \bar\eta_{n-2})\in{\cal P}_\frg$ implies the following equalities: \begin{eqnarray*}\partial(\xi_1\lrcorner\dots\lrcorner\xi_{n-2}\lrcorner\bar\eta_1\lrcorner\dots\lrcorner\bar\eta_{n-2}\lrcorner dV) = 0  \hspace{2ex} \mbox{and} \hspace{2ex} \bar\partial(\xi_1\lrcorner\dots\lrcorner\xi_{n-2}\lrcorner\bar\eta_1\lrcorner\dots\lrcorner\bar\eta_{n-2}\lrcorner dV) = 0.\end{eqnarray*} Thanks to (\ref{eqn:moment-map_w-d_proof_2}), this shows that \begin{eqnarray*}\int_S(\partial\bar\beta + \bar\partial\beta)(\xi_1,\dots , \xi_{n-2},\,\bar\eta_1,\dots , \bar\eta_{n-2})\,dV = 0\end{eqnarray*} for every $(\xi_1,\dots , \xi_{n-2},\,\bar\eta_1,\dots , \bar\eta_{n-2})\in{\cal P}_\frg$. This proves (\ref{eqn:moment-map_w-d_proof_1}) and completes the proof of Proposition \ref{Prop:moment-map_well-defined}.  \hfill $\Box$

\section{Removal of the assumption $H^{n-2,\,n-2}_A(S,\,\C)=\{0\}$}\label{section:removal}

  In Proposition \ref{Prop:moment-map_well-defined}, the assumption $H^{n-2,\,n-2}_A(S,\,\C)=\{0\}$ was made in order to ensure that the object introduced in Definition \ref{Def:moment-map_def} is independent of the choice of the solution $\Gamma_f$ of equation (\ref{eqn:f-omega_n1_potential}). This solution is unique only up to $\ker(\partial\bar\partial)$, but it can be made unique in the absolute sense by choosing it to have minimal $L^2$-norm with respect to a given Hermitian metric on $S$.

  To avoid fixing an arbitrary Hermitian metric on $S$ that would be unrelated to the already fixed volume form $dV>0$, we will use the following generalisation to the Hermitian case, given by Tosatti and Weinkove in [TW10] as the culmination of a string of works, including [Che87] and [GL09], of Yau's celebrated resolution [Yau78] of the Calabi Conjecture for the K\"ahler case.

\begin{The}([TW10])\label{The:TW} Let $S$ be a compact complex manifold with $\mbox{dim}_{\C}S=d$ and let $\gamma$ be a Hermitian metric on $S$. 

 Then, for any $C^{\infty}$ function $F\,:\,S\rightarrow \R$, there exist a unique constant $C>0$ and a unique $C^{\infty}$ function $\varphi\,:\,S\rightarrow \R$ such that \begin{eqnarray}\label{eqn:TW_M-A_eqn}(\gamma + i\partial\bar\partial\varphi)^d = C e^F\gamma^d,\hspace{2ex} \gamma + i\partial\bar\partial\varphi >0 \hspace{2ex} \mbox{and} \hspace{2ex} \sup\limits_S\varphi=0.\end{eqnarray}

\end{The}

\vspace{2ex}

For every $C^\infty$ $(1,\,1)$-form $\alpha$ on $S$, we consider the following set: \begin{eqnarray*}[\alpha]_{BC}:=\bigg\{\alpha + i\partial\bar\partial\varphi\,\mid\,\varphi:S\longrightarrow\R \hspace{2ex} C^\infty\hspace{0.5ex}\mbox{function}\bigg\}\end{eqnarray*} that we call the {\it pseudo-Bott-Chern class} of $\alpha$. If $d\alpha=0$, $[\alpha]_{BC}$ is the genuine Bott-Chern cohomology class of $\alpha$.

  Note that, unlike the Bott-Chern cohomology group $H^{1,\,1}_{BC}(S,\,\C)$ (consisting of all Bott-Chern classes) of $S$, the $\C$-vector space $\widetilde{H^{1,\,1}_{BC}}(S,\,\C)$ consisting of all the {\it pseudo-Bott-Chern classes} of $S$ is {\it infinite-dimensional}. To see this, fix an arbitrary Hermitian metric $\gamma_0$ on $S$ and consider the standard $3$-space $L^2_{\gamma_0}$-orthogonal decomposition: \begin{eqnarray}\label{eqn:3-space_B-C}C^\infty_{1,\,1}(S,\,\C) = {\cal H}^{1,\,1}_{\gamma_0}(S,\,\C)\oplus\mbox{Im}\,(\partial\bar\partial)\oplus\bigg(\mbox{Im}\,\partial^\star + \mbox{Im}\,\bar\partial^\star\bigg),\end{eqnarray} where ${\cal H}^{1,\,1}_{\gamma_0}(S,\,\C)$ is the Bott-Chern-harmonic space of bidegree $(1,\,1)$ induced by $\gamma_0$, namely the kernel of the Bott-Chern Laplacian \begin{equation*}\Delta_{BC}=\partial^{\star}\partial + \bar\partial^{\star}\bar\partial + (\partial\bar\partial)^{\star}(\partial\bar\partial) + (\partial\bar\partial)(\partial\bar\partial)^{\star} + (\partial^{\star}\bar\partial)^{\star}(\partial^{\star}\bar\partial) + (\partial^{\star}\bar\partial)(\partial^{\star}\bar\partial)^{\star}:C^\infty_{1,\,1}(S,\,\C)\longrightarrow C^\infty_{1,\,1}(S,\,\C)\end{equation*} and $\partial^{\star} =\partial^\star_{\gamma_0}$, $\bar\partial^{\star} =\bar\partial^\star_{\gamma_0}$ are the adjoints of $\partial$ and $\bar\partial$ with respect to the $L^2$-inner product defined by $\gamma_0$. As is well known by standard harmonic theory, the ellipticity of $\Delta_{BC}$ and the compactness of $S$ imply the finite-dimensionality of ${\cal H}^{1,\,1}_{\gamma_0}(S,\,\C)$. However, the vector spaces $\mbox{Im}\,\partial^\star$ and $\mbox{Im}\,\bar\partial^\star$ are infinite-dimensional. Meanwhile, for every $\alpha\in C^\infty_{1,\,1}(S,\,\C)$, if \begin{eqnarray*}\alpha = \alpha_h + i\partial\bar\partial\varphi_\alpha + \bigg(\partial^\star\beta_\alpha^{2,\,1} + \bar\partial^\star\beta_\alpha^{1,\,2}\bigg)\end{eqnarray*} is the splitting of $\alpha$ according to (\ref{eqn:3-space_B-C}), we get: \begin{eqnarray*}[\alpha]_{BC} = [\alpha_h]_{BC} + [\partial^\star\beta_\alpha^{2,\,1} + \bar\partial^\star\beta_\alpha^{1,\,2}]_{BC}.\end{eqnarray*}

  Now, when $\alpha$ varies across $C^\infty_{1,\,1}(S,\,\C)$, the Bott-Chern class $[\alpha_h]_{BC}$ varies across the finite-dimensional Bott-Chern cohomology space $H^{1,\,1}(S,\,\C)$ of $S$, but the pseudo-Bott-Chern class $[\partial^\star\beta_\alpha^{2,\,1} + \bar\partial^\star\beta_\alpha^{1,\,2}]_{BC}$ varies across the infinite-dimensional vector space $\mbox{Im}\,\partial^\star + \mbox{Im}\,\bar\partial^\star$ (since $\mbox{Im}\,(\partial\bar\partial)$ meets $\mbox{Im}\,\partial^\star + \mbox{Im}\,\bar\partial^\star$ only at zero).

  Besides the finite-dimensional vector subspace $H^{1,\,1}_{BC}(S,\,\C)$, the infinite-dimensional vector space $\widetilde{H^{1,\,1}_{BC}}(S,\,\C)$ contains the cone \begin{eqnarray*}{\cal M}(S):=\bigg\{[\gamma]_{BC}\in\widetilde{H^{1,\,1}_{BC}}(S,\,\C)\,\mid\,\gamma>0\bigg\}\end{eqnarray*} consisting of all the pseudo-Bott-Chern classes that admit a {\it positive definite} representative $\gamma$.

\vspace{2ex}

Now, in our situation, since a volume form $dV>0$ has been fixed on $S$, for every pseudo-Bott-Chern class $\mathfrak{c} = [\gamma]_{BC}\in{\cal M}(S)$ representable by a Hermitian metric $\gamma$ on $S$, the Tosatti-Weinkove Theorem \ref{The:TW} ensures the existence of a unique constant $C>0$ and of a unique $C^\infty$ function $\varphi:S\longrightarrow\R$ such that \begin{eqnarray}\label{eqn:TW_M-A_eqn_application}(\gamma + i\partial\bar\partial\varphi)^d = C\,dV,\hspace{2ex} \gamma + i\partial\bar\partial\varphi >0 \hspace{2ex} \mbox{and} \hspace{2ex} \sup\limits_S\varphi=0.\end{eqnarray}

In other words, $\gamma_\varphi:= \gamma + i\partial\bar\partial\varphi$ is the unique Hermitian metric lying in the given pseudo-Bott-Chern class $\mathfrak{c} = [\gamma]_{BC}\in{\cal M}(S)$ whose volume form $\gamma_\varphi^d$ is a constant multiple of $dV$.

On the other hand, the Hermitian metric $\gamma_\varphi$ induces its Bott-Chern Laplacian 
$$\Delta_{BC}: C^\infty_{n-1,\,n-1}(S,\,\C) \longrightarrow C^\infty_{n-1,\,n-1}(S,\,\C)$$ 
defined by the formula recalled above, while Theorem 4.1 in [Pop15] gives the following Neumann-type formula for the (unique) solution $\Gamma_f$ of the $\partial\bar\partial$-equation (\ref{eqn:f-omega_n1_potential}): \begin{equation}\label{eqn:Neumann_ddbar}i\Gamma_f = (\partial\bar\partial)^\star\Delta_{BC}^{-1}(f^\star\omega_{n-1}),\end{equation} where $\Delta_{BC}^{-1}$ is the Green operator of $\Delta_{BC}$. Since this form $\Gamma_f$ depends on the given pseudo-Bott-Chern class $\mathfrak{c} = [\gamma]_{BC}\in{\cal M}(S)$ and is uniquely induced by $f$ and $\mathfrak{c}$, we will denote it by $\Gamma_{f,\,\mathfrak{c}}$.

\vspace{2ex}

In view of this discussion, we can adapt Definition \ref{Def:moment-map_def} in the following way.

\begin{Def}\label{Def:moment-map_BC_def} With every $n$-dimensional complex {\bf balanced} manifold $(X,\,\omega_{n-1})$ and every $d$-dimensional {\bf compact} complex manifold $(S,\,dV,\,\mathfrak{c})$ equipped with a positive {\bf volume form} $dV>0$ and a pseudo-Bott-Chern class $\mathfrak{c}\in{\cal M}(S)$ representable by a Hermitian metric such that $d\geq n-1$, we associate the map $\mu_{\mathfrak{c}}:{\mathscr X}\longrightarrow{\cal P}_\frg^\star$ defined by \begin{eqnarray}\label{eqn:moment-map_BC_def}\mu_{\mathfrak{c}}(f)(\xi_1,\dots , \xi_{n-2},\,\bar\eta_1,\dots , \bar\eta_{n-2}):=i\int\limits_S\Gamma_{f,\,\mathfrak{c}}(\xi_1,\dots , \xi_{n-2},\,\bar\eta_1,\dots , \bar\eta_{n-2})\,dV\end{eqnarray} for every $f\in{\mathscr X}$ and every $(\xi_1,\dots , \xi_{n-2},\,\bar\eta_1,\dots , \bar\eta_{n-2})\in{\cal P}_\frg$.

\end{Def}

For this definition, no analogue of Proposition \ref{Prop:moment-map_well-defined} is necessary since the form $\Gamma_{f,\,\mathfrak{c}}$ is uniquely determined by $f$ and $\mathfrak{c}$.

\section{Lie derivatives of types $(1,\,0)$ and $(0,\,1)$}\label{section:Lie-derivatives} We take this opportunity to introduce the complex analogues in bidegrees $(1,\,0)$ and $(0,\,1)$ of the Lie derivative w.r.t. a vector field of the real case. The following definition was obliquely suggested by Sarkaria in [Sar78, $\S.7$].

\begin{Def}\label{Def:Lie-deriv_def} Let $X$ be a complex manifold.

\vspace{1ex}  

(i)\, For any smooth vector field $\xi\in C^{\infty}(X,\,T^{1,\,0}X)$ of type $(1,\,0)$, the {\bf $(1,\,0)$-Lie derivative} w.r.t. $\xi$ is defined as

\begin{equation}\label{eqn:Lie-deriv_1-0_def}L^{1,\,0}_{\xi}:=[\xi\lrcorner\cdot,\,\partial].\end{equation}

\noindent This means that on any differential form $u$, we have $L^{1,\,0}_{\xi}u:= \xi\lrcorner(\partial u) + \partial(\xi\lrcorner u)$.

\vspace{1ex}  

(ii)\, For any smooth vector field $\eta\in C^{\infty}(X,\,T^{1,\,0}X)$ of type $(1,\,0)$, the {\bf $(0,\,1)$-Lie derivative} w.r.t. $\bar\eta$ is defined as

\begin{equation}\label{eqn:Lie-deriv_0-1_def}L^{0,\,1}_{\bar\eta}:=[\bar\eta\lrcorner\cdot,\,\bar\partial].\end{equation}

\noindent This means that on any differential form $u$, we have $L^{0,\,1}_{\bar\eta}u:= \bar\eta\lrcorner(\bar\partial u) + \bar\partial(\bar\eta\lrcorner u)$.

\end{Def}

Here, as throughout the text, we use the standard notation $[\,\,,\,\,]$ for the graded commutator of two endomorphisms $A, B$ of respective degrees $a,b$ of the graded algebra $C^{\infty}_{\bullet}(X,\,\C)$ of smooth differential forms on $X$: $[A,\,B]=AB-(-1)^{ab}BA$.

Note that, since $T^{1,\,0}X$ is a {\it holomorphic} vector bundle, it has a canonical $\bar\partial$-operator associated with its holomorphic structure, so $\bar\partial\xi\in C^\infty_{0,\,1}(X,\,T^{1,\,0}X)$ is a well-defined smooth $(0,\,1)$-form with values in $T^{1,\,0}X$ for every $\xi\in C^{\infty}(X,\,T^{1,\,0}X)$. In particular, there is no need for a Lie derivative of type $(0,\,1)$ w.r.t. $\xi$ since, as can be immediately checked in coordinates, the expected Leibniz formula holds: \begin{eqnarray*}\bar\partial(\xi\lrcorner u) = \bar\partial\xi\lrcorner u - \xi\lrcorner\bar\partial u,  \hspace{3ex} \xi\in C^{\infty}(X,\,T^{1,\,0}X),\end{eqnarray*} for any scalar-valued differential form $u$ on $X$, where $\bar\partial\xi\lrcorner u$ is defined by contracting $u$ with the vector-field part of $\bar\partial\xi$ and multiplying it by the differential-form part of $\bar\partial\xi$. For example, if $u$ is of type $(p,\,q)$, $\bar\partial\xi\lrcorner u$ is a scalar-valued form of type $(p-1,\,q+1)$. By contrast, the holomorphic vector bundle $T^{1,\,0}X$ has no canonical $\partial$-operator (i.e. no canonical $(1,\,0)$-connection), so $\partial\xi$ is meaningless. The $(1,\,0)$-Lie derivative $L^{1,\,0}_\xi$ defined above plays the role of what would have been the contraction-multiplication by $\partial\xi$ to produce the Leibniz-type formula of (i) of Definition \ref{Def:Lie-deriv_def}.

Symmetrically, $T^{0,\,1}X$ is an {\it anti-holomorphic} vector bundle, so it has a canonical $\partial$-operator (and, in particular, $\partial\bar\eta\in C^\infty_{1,\,0}(X,\,T^{0,\,1}X)$ is a well-defined smooth $(1,\,0)$-form with values in $T^{0,\,1}X$ for every smooth $(0,\,1)$-vector field $\bar\eta$ on $X$), but it has no canonical $\bar\partial$-operator. It can be immediately checked in coordinates that the expected Leibniz formula holds: \begin{eqnarray*}\partial(\bar\eta\lrcorner u) = \partial\bar\eta\lrcorner u - \bar\eta\lrcorner\partial u,  \hspace{3ex} \bar\eta\in C^{\infty}(X,\,T^{0,\,1}X),\end{eqnarray*} for any scalar-valued differential form $u$ on $X$. Meanwhile, the $(0,\,1)$-Lie derivative $L^{0,\,1}_{\bar\eta}$ plays the role of what would have been the contraction-multiplication by the non-existent $\bar\partial\bar\eta$.

\vspace{1ex}

We now notice that $L^{1,\,0}_{\xi}$ and $L^{0,\,1}_{\bar\eta}$ have analogous properties to those of the standard Lie derivative of the real case. We first deal with the case of $L^{1,\,0}_{\xi}$.

\begin{Lem}\label{Lem:Lie-deriv-prop_1-0} Let $\xi,\eta\in C^{\infty}(X,\,T^{1,\,0}X)$.

\vspace{1ex}

  $(i)$\, The standard Lie derivative $L_{\xi}= [\xi\lrcorner\cdot,\,d] $ w.r.t. $\xi$ is related to $L^{1,\,0}_{\xi}$ as follows: \begin{equation}\label{eqn:Lie-deriv-relation_1-0}L_{\xi}=L^{1,\,0}_{\xi} + (\bar\partial\xi)\lrcorner\cdot.\end{equation}

\noindent In particular, if $\xi$ is holomorphic, then $L_{\xi}=L^{1,\,0}_{\xi}$.

\noindent $(ii)$\, For any function $f\in C^{\infty}(X,\,\C)$, we have $L^{1,\,0}_{\xi}f = \xi\cdot f$.

\vspace{1ex}

\noindent $(iii)$\, The following identities hold: \begin{equation*}(a)\,\, [L^{1,\,0}_{\xi},\,\partial]=0 \hspace{2ex} \mbox{and} \hspace{2ex} (b)\,\, [L^{1,\,0}_{\xi},\,\bar\partial]= [\partial,\,\bar\partial\xi\lrcorner\cdot].\end{equation*} In particular, if $\xi$ is holomorphic, then $[L^{1,\,0}_{\xi},\,\bar\partial] = 0$.

\vspace{1ex}

\noindent $(iv)$\, The following identities hold: \begin{equation}\label{eqn:L_commutations}(a)\,\,[\xi\lrcorner\cdot,\,L^{1,\,0}_{\eta}] = [L^{1,\,0}_{\xi},\,\eta\lrcorner\cdot] = [\xi,\,\eta]\lrcorner\cdot \hspace{2ex} \mbox{and} \hspace{2ex} (b)\,\,[L^{1,\,0}_{\xi},\,L^{1,\,0}_{\eta}] = L^{1,\,0}_{[\xi,\,\eta]}.\end{equation}

\noindent $(v)$\, For any differential forms $u,v$ (of any degrees), we have $L^{1,\,0}_{\xi}(u\wedge v) = (L^{1,\,0}_{\xi}u)\wedge v + u\wedge L^{1,\,0}_{\xi}v$.

\noindent $(vi)$\, If $\xi$ is holomorphic, then \begin{equation*}L^{1,\,0}_\xi(\bar\eta\lrcorner\alpha) = \bar\eta\lrcorner L^{1,\,0}_\xi(\alpha) + [\xi,\,\bar\eta]\lrcorner\alpha\end{equation*} for any differential form $\alpha$.

\noindent $(vii)$\, For any smooth functions $f, g$ and any smooth $(0,\,q)$-form $\alpha$ on $X$, we have: $$L^{1,\,0}_{f\xi + g\eta}(\alpha) = fL^{1,\,0}_\xi(\alpha) + gL^{1,\,0}_\eta(\alpha).$$

\end{Lem}

\noindent {\it Proof.} $(i)$\, Let $u$ be any smooth differential form on $X$ (of any degree). Then $\partial(\xi\lrcorner u) = -\xi\lrcorner\partial u + L^{1,\,0}_{\xi}u$ and $\bar\partial(\xi\lrcorner u) = -\xi\lrcorner\bar\partial u + (\bar\partial\xi)\lrcorner u$. Since $d(\xi\lrcorner u) = -\xi\lrcorner du + L_{\xi}u$, (\ref{eqn:Lie-deriv-relation_1-0}) follows.

\hspace{1ex}

$(ii)$\, $L^{1,\,0}_{\xi}f=\xi\lrcorner\partial f + \partial(\xi\lrcorner f) = \xi\lrcorner\partial f$ since $\xi\lrcorner f=0$ as a $(-1,\,0)$-form. If, in local coordinates, $\xi=\sum\xi_j\frac{\partial}{\partial z_j}$, then $\xi\lrcorner\partial f = \sum\xi_j\frac{\partial f}{\partial z_j} = \xi\cdot f$.

\hspace{1ex}

$(iii)$\,(a)\, To process the expression $[L^{1,\,0}_{\xi},\,\partial] = [[\xi\lrcorner\cdot,\,\partial],\,\partial]$, we use the Jacobi identity \begin{eqnarray*}[[\xi\lrcorner\cdot,\,\partial],\,\partial] + [[\partial,\,\partial],\,\xi\lrcorner\cdot] + [[\partial,\,\xi\lrcorner\cdot],\,\partial] = 0.\end{eqnarray*} Since $[\partial,\,\partial]=0$ and $[\partial,\,\xi\lrcorner\cdot] = [\xi\lrcorner\cdot,\,\partial]$, this translates to $2[[\xi\lrcorner\cdot,\,\partial],\,\partial]=0$, hence to $[L^{1,\,0}_{\xi},\,\partial] = 0$. 

\vspace{1ex}

 (b)\, To process the expression $[L^{1,\,0}_{\xi},\,\bar\partial] = [[\xi\lrcorner\cdot,\,\partial],\,\bar\partial]$, we use the Jacobi identity \begin{eqnarray*}[[\xi\lrcorner\cdot,\,\partial],\,\bar\partial] + [[\partial,\,\bar\partial],\,\xi\lrcorner\cdot] + [[\bar\partial,\,\xi\lrcorner\cdot],\,\partial] = 0.\end{eqnarray*} This proves the contention since $[\partial,\,\bar\partial]= \partial\bar\partial + \bar\partial\partial = 0$ and $[\bar\partial,\,\xi\lrcorner\cdot] = (\bar\partial\xi)\lrcorner\cdot$, while $[(\bar\partial\xi)\lrcorner\cdot,\,\partial] = - [\partial,\,(\bar\partial\xi)\lrcorner\cdot]$.

\vspace{1ex}

 $(iv)$\, By definition of $L^{1,\,0}_{\eta}$, we have $[\xi\lrcorner\cdot,\,L^{1,\,0}_{\eta}] = [\xi\lrcorner\cdot,\,[\eta\lrcorner\cdot,\,\partial]]$. The Jacobi identity gives

$$-[\xi\lrcorner\cdot,\,[\eta\lrcorner\cdot,\,\partial]] - [\eta\lrcorner\cdot,\,[\partial,\,\xi\lrcorner\cdot]] - [\partial,\,[\xi\lrcorner\cdot,\,\eta\lrcorner\cdot]] =0.$$

\noindent Since $[\xi\lrcorner\cdot,\,\eta\lrcorner\cdot]=0$, we get $[\xi\lrcorner\cdot,\,L^{1,\,0}_{\eta}] = - [\eta\lrcorner\cdot,\,L^{1,\,0}_{\xi}] = [L^{1,\,0}_{\xi},\,\eta\lrcorner\cdot]$. This proves the first identity in $(a)$. The identity $[\xi\lrcorner\cdot,\,L^{1,\,0}_{\eta}] = [\xi,\,\eta]\lrcorner\cdot$ is equivalent to the following well-known intrinsic formula (the so-called Cartan formula) for the exterior derivative $\partial$. For any $\alpha\in C^{\infty}(X,\,\Lambda^{k,\,0}T^{\star}X)$ and any $\xi_0, \dots , \xi_k\in C^{\infty}(X,\,T^{1,\,0}X)$, we have: \begin{eqnarray}\label{eqn:del-intrinsic}\nonumber (\partial\alpha)(\xi_0,\dots , \xi_k) & = & \sum\limits_{j=0}^k(-1)^j\,\xi_j\cdot\alpha(\xi_0,\dots , \widehat{\xi_j},\dots , \xi_k) \\
      & + & \sum\limits_{0\leq j<l\leq k}(-1)^{j+l}\,\alpha([\xi_j,\,\xi_l],\,\xi_0,\dots , \widehat{\xi_j},\dots , \widehat{\xi_l},\dots , \xi_k).\end{eqnarray} 

\noindent Indeed, let us check this equivalence for $(1,\,0)$-forms $\alpha$. The identity $[\xi\lrcorner\cdot,\,L^{1,\,0}_{\eta}]\,\alpha = [\xi,\,\eta]\lrcorner\alpha$ is equivalent to $-\xi\lrcorner L^{1,\,0}_{\eta}\alpha = - \eta\cdot(\xi\lrcorner\alpha) - [\xi,\,\eta]\lrcorner\alpha$ (see $(ii)$), i.e. to $-\xi\lrcorner[\eta\lrcorner\partial\alpha + \partial(\eta\lrcorner\alpha)] = - \eta\cdot\alpha(\xi) - \alpha([\xi,\,\eta])$, which in turn is equivalent to $(\partial\alpha)(\xi,\,\eta) - \xi\cdot\alpha(\eta) = - \eta\cdot\alpha(\xi) - \alpha([\xi,\,\eta])$. This completes the proof of $(a)$ in (\ref{eqn:L_commutations}). 

 To prove $(b)$, notice that $[L^{1,\,0}_{\xi},\,L^{1,\,0}_{\eta}] = [[\xi\lrcorner\cdot,\,\partial],\,L^{1,\,0}_{\eta}]$ and the Jacobi identity spells

$$[[\xi\lrcorner\cdot,\,\partial],\,L^{1,\,0}_{\eta}] - [[\partial,\,L^{1,\,0}_{\eta}],\,\xi\lrcorner\cdot] + [[L^{1,\,0}_{\eta},\,\xi\lrcorner\cdot],\,\partial] = 0.$$ 

 \noindent Since $[\partial,\,L^{1,\,0}_{\eta}]=0$ by $(iii)$ and $[L^{1,\,0}_{\eta},\,\xi\lrcorner\cdot] = [\eta,\,\xi]\lrcorner\cdot = - [\xi,\,\eta]\lrcorner\cdot$ by $(a)$, the Jacobi identity translates to $[L^{1,\,0}_{\xi},\,L^{1,\,0}_{\eta}] = [[\xi,\,\eta]\lrcorner\cdot,\,\partial]$ which is precisely $(b)$ of (\ref{eqn:L_commutations}).

 \hspace{1ex}

 $(v)$\, From the definition of $L^{1,\,0}_{\xi}$, we get

\noindent $L^{1,\,0}_{\xi}(u\wedge v) = \xi\lrcorner\partial(u\wedge v) + \partial(\xi\lrcorner(u\wedge v)) = \xi\lrcorner[\partial u\wedge v + (-1)^{deg\,u}u\wedge\partial v] + \partial[(\xi\lrcorner u)\wedge v + (-1)^{deg\,u}u\wedge(\xi\lrcorner v)]$

\vspace{1ex}

\hspace{7ex} $ = (\xi\lrcorner\partial u)\wedge v + \partial(\xi\lrcorner u)\wedge v + u\wedge(\xi\lrcorner\partial v) + u\wedge\partial(\xi\lrcorner v) = (L^{1,\,0}_{\xi}u)\wedge v + u\wedge L^{1,\,0}_{\xi}v.$

\hspace{1ex}

$(vi)$\, Using a standard property of the usual Lie derivative $L_\xi$ w.r.t. $\xi$, we get: \begin{eqnarray*}L_\xi(\bar\eta\lrcorner\alpha) = \bar\eta\lrcorner L_\xi(\alpha) + [\xi,\,\bar\eta]\lrcorner\alpha.\end{eqnarray*} Using property (i), this translates to \begin{eqnarray*}L^{1,\,0}_\xi(\bar\eta\lrcorner\alpha) + (\bar\partial\xi)\lrcorner\bar\eta\lrcorner\alpha = \bar\eta\lrcorner L^{1,\,0}_\xi(\alpha) + \bar\eta\lrcorner(\bar\partial\xi)\lrcorner\alpha + [\xi,\,\bar\eta]\lrcorner\alpha.\end{eqnarray*} In our case, $\bar\partial\xi=0$ since $\xi$ is assumed holomorphic, so we get the contention.    

\hspace{1ex}

$(vii)$\, The additivity property $L^{1,\,0}_{f\xi + g\eta}(\alpha) = L^{1,\,0}_{f\xi}(\alpha) + L^{1,\,0}_{g\eta}(\alpha)$ being obvious, from the definition of $L^{1,\,0}$, we will only prove the homogeneity property $L^{1,\,0}_{f\xi}(\alpha) = fL^{1,\,0}_{\xi}(\alpha)$.   

The definition of $L^{1,\,0}_\xi$ gives: $\partial(\xi\lrcorner\alpha) = L^{1,\,0}_\xi(\alpha) - \xi\lrcorner\partial\alpha$. Hence, \begin{eqnarray*}f\partial(\xi\lrcorner\alpha) = fL^{1,\,0}_\xi(\alpha) - f(\xi\lrcorner\partial\alpha).\end{eqnarray*}

Meanwhile, the definition of $L^{1,\,0}_{f\xi}$ gives the first equality in the equivalence sequence below: \begin{eqnarray*} & & \partial((f\xi)\lrcorner\alpha) = L^{1,\,0}_{f\xi}(\alpha) - (f\xi)\lrcorner\partial\alpha \iff f\partial(\xi\lrcorner\alpha) + (\partial f)\wedge(\xi\lrcorner\alpha) = L^{1,\,0}_{f\xi}(\alpha) - f(\xi\lrcorner\partial\alpha) \\
 & \iff & fL^{1,\,0}_\xi(\alpha) - f(\xi\lrcorner\partial\alpha) + (\partial f)\wedge(\xi\lrcorner\alpha) = L^{1,\,0}_{f\xi}(\alpha) - f(\xi\lrcorner\partial\alpha)   \end{eqnarray*} This proves the expected equality $L^{1,\,0}_{f\xi}(\alpha) = fL^{1,\,0}_{\xi}(\alpha)$ after $f(\xi\lrcorner\partial\alpha)$ drops out from both sides and after we take into account that $\xi\lrcorner\alpha=0$ since $\alpha$ is of type $(0,\,q)$.  \hfill $\Box$

\vspace{3ex}

The analogous properties of $L^{0,\,1}_{\bar\eta}$ are obtained by conjugation from those of $L^{1,\,0}_{\xi}$ listed in the previous lemma.   

\begin{Lem}\label{Lem:Lie-deriv-prop_0-1} Let $\xi,\eta\in C^{\infty}(X,\,T^{1,\,0}X)$.

\vspace{1ex}

  $(i)$\, The standard Lie derivative $L_{\bar\eta}= [\bar\eta\lrcorner\cdot,\,d]$ w.r.t. $\bar\eta$ is related to $L^{0,\,1}_{\bar\eta}$ as follows: \begin{equation}\label{eqn:Lie-deriv-relation_0-1}L_{\bar\eta}=L^{0,\,1}_{\bar\eta} + (\partial\bar\eta)\lrcorner\cdot.\end{equation}

\noindent In particular, if $\eta$ is holomorphic, then $L_{\bar\eta}=L^{0,\,1}_{\bar\eta}$.

\noindent $(ii)$\, For any function $f\in C^{\infty}(X,\,\C)$, we have $L^{0,\,1}_{\bar\eta}f = \bar\eta\cdot f$.

\vspace{1ex}

\noindent $(iii)$\, The following identities hold: \begin{equation*}(a)\,\, [L^{0,\,1}_{\bar\eta},\,\bar\partial]=0 \hspace{2ex} \mbox{and} \hspace{2ex} (b)\,\, [L^{0,\,1}_{\bar\eta},\,\partial]= [\bar\partial,\,\partial\bar\eta\lrcorner\cdot].\end{equation*} In particular, if $\eta$ is holomorphic, then $[L^{0,\,1}_{\bar\eta},\,\partial] = 0$.

\vspace{1ex}

\noindent $(iv)$\, The following identities hold: \begin{equation}\label{eqn:L_commutations-1}(a)\,\,[\bar\xi\lrcorner\cdot,\,L^{0,\,1}_{\bar\eta}] = [L^{0,\,1}_{\bar\xi},\,\bar\eta\lrcorner\cdot] = [\bar\xi,\,\bar\eta]\lrcorner\cdot \hspace{2ex} \mbox{and} \hspace{2ex} (b)\,\,[L^{0,\,1}_{\bar\xi},\,L^{0,\,1}_{\bar\eta}] = L^{0,\,1}_{[\bar\xi,\,\bar\eta]}.\end{equation}

\noindent $(v)$\, For any differential forms $u,v$ (of any degrees), we have $L^{0,\,1}_{\bar\eta}(u\wedge v) = (L^{0,\,1}_{\bar\eta}u)\wedge v + u\wedge L^{0,\,1}_{\bar\eta}v$.

\noindent $(vi)$\, If $\eta$ is holomorphic, then \begin{equation*}L^{0,\,1}_{\bar\eta}(\xi\lrcorner\alpha) = \xi\lrcorner L^{0,\,1}_{\bar\eta}(\alpha) + [\bar\eta,\,\xi]\lrcorner\alpha\end{equation*} for any differential form $\alpha$.

\noindent $(vii)$\, For any smooth functions $f, g$ and any smooth $(p,\,0)$-form $\alpha$ on $X$, we have: $$L^{0,\,1}_{f\bar\xi + g\bar\eta}(\alpha) = fL^{0,\,1}_{\bar\xi}(\alpha) + gL^{0,\,1}_{\bar\eta}(\alpha).$$

\end{Lem}

\section{Moment map confirmation}\label{section:moment-map_confirmed} The set-up is the one described in $\S.$\ref{section:moment-map_w-d} and $\S.$\ref{section:removal}. For the sake of unifying the notation, we will still denote by $\mu$ the map $\mu_{\mathfrak{c}}$ of Definition \ref{Def:moment-map_BC_def} and by $\Gamma_f$ the form $\Gamma_{f,\,\mathfrak{c}}$. 

In the context of either Definition \ref{Def:moment-map_def} or Definition \ref{Def:moment-map_BC_def}, for every $(\xi_1,\dots , \xi_{n-2},\,\bar\eta_1,\dots , \bar\eta_{n-2})\in{\cal P}_\frg$, we consider the scalar-valued map:  \begin{eqnarray*}\label{eqn:moment-map_def_confirmed}\nonumber\langle\mu(\cdot),\,(\xi_1,\dots , \xi_{n-2},\,\bar\eta_1,\dots , \bar\eta_{n-2})\rangle:{\mathscr X} \longrightarrow  \C,\end{eqnarray*}
\begin{eqnarray}\label{eqn:moment-map_def_confirmed}  f \longmapsto \langle\mu(f),\,(\xi_1,\dots , \xi_{n-2},\,\bar\eta_1,\dots , \bar\eta_{n-2})\rangle:=i\int\limits_S\Gamma_f(\xi_1,\dots , \xi_{n-2},\,\bar\eta_1,\dots , \bar\eta_{n-2})\,dV.\end{eqnarray} Thus, $\mu(f)(\xi_1,\dots , \xi_{n-2},\,\bar\eta_1,\dots , \bar\eta_{n-2}) = \langle\mu(f),\,(\xi_1,\dots , \xi_{n-2},\,\bar\eta_1,\dots , \bar\eta_{n-2})\rangle$.

\begin{The}\label{The:moment-map_def_final} Let $(X,\,\omega_{n-1})$ be an $n$-dimensional complex {\bf balanced} manifold and let $(S,\,dV)$ be a $d$-dimensional {\bf compact} complex manifold equipped with a positive {\bf volume form} such that $d\geq n-1$. Suppose that

\vspace{1ex}

$\bullet$ either $H^{n-2,\,n-2}_A(S,\,\C)=\{0\}$;

\vspace{1ex}

$\bullet$ or $\mathfrak{c}\in{\cal M}(S)$ is a pseudo-Bott-Chern class representable by a Hermitian metric on $S$.

\vspace{1ex}

   Then, the map $\mu:{\mathscr X}\longrightarrow{\cal P}_\frg^\star$ introduced in either Definition \ref{Def:moment-map_def} or Definition \ref{Def:moment-map_BC_def} is a {\bf moment map} for the action of ${\cal G}$ on $({\mathscr X},\,\Omega)$ in the following sense: for every $(\xi_1,\dots , \xi_{n-2},\,\bar\eta_1,\dots , \bar\eta_{n-2})\in{\cal P}_\frg$, \begin{eqnarray}\label{eqn:moment-map_def_final}\partial\bar\partial\bigg(\langle\mu(\cdot),\,(\xi_1,\dots , \xi_{n-2},\,\bar\eta_1,\dots , \bar\eta_{n-2})\rangle\bigg) = X(\bar\eta_{n-2})\lrcorner\dots\lrcorner X(\bar\eta_1)\lrcorner X(\xi_{n-2})\lrcorner\dots\lrcorner X(\xi_1)\lrcorner\Omega,\end{eqnarray} where $X(\xi_j)$, resp. $X(\bar\eta_k)$, is the vector field on ${\mathscr X}$ defined by the infinitesimal action of $\xi_j$, resp. $\bar\eta_k$.

\end{The}

Before giving the proof, we make a trivial observation. Let $\omega$ be a $(1,\,1)$-form on a complex manifold. In local coordinates, it reads: $\omega = \sum_{j,\,k}\omega_{j\bar{k}}\,idz_j\wedge d\bar{z}_k$, so $\omega(\partial/\partial z_l,\,\partial/\partial\bar{z}_r) = i\omega_{l\bar{r}} = (\partial/\partial\bar{z}_r)\lrcorner(\partial/\partial z_l)\lrcorner\omega$. Meanwhile, if $\omega = i\partial\bar\partial f$, then  \begin{eqnarray*}\omega_{l\bar{r}} = \frac{\partial^2f}{\partial z_l\,\partial\bar{z}_r} = L^{1,\,0}_{\frac{\partial}{\partial z_l}}L^{0,\,1}_{\frac{\partial}{\partial\bar{z}_r}}(f),  \hspace{6ex} l,r,\end{eqnarray*} hence \begin{eqnarray}\label{eqn:iddbar-f_2-contraction_Lie_coord}(\partial/\partial\bar{z}_r)\lrcorner(\partial/\partial z_l)\lrcorner(i\partial\bar\partial f) = (i\partial\bar\partial f)(\partial/\partial z_l,\,\partial/\partial\bar{z}_r) = iL^{1,\,0}_{\frac{\partial}{\partial z_l}}L^{0,\,1}_{\frac{\partial}{\partial\bar{z}_r}}(f),  \hspace{6ex} l,r.\end{eqnarray}

The following observation extends this to general vector fields. 

\begin{Lem}\label{Lem:eqn:iddbar-f_2-contraction_Lie} Let $f$ be a smooth function on a complex manifold $X$ and let $v,w$ be smooth vector fields of type $(1,\,0)$ on $X$. 

\vspace{1ex}

(i)\, If $w$ is {\bf holomorphic}, then $(i\partial\bar\partial f)(v,\,\bar{w}) = iL^{1,\,0}_vL^{0,\,1}_{\bar{w}}(f)$.

\vspace{1ex}

(ii)\, If $v$ is {\bf holomorphic}, then $(i\partial\bar\partial f)(v,\,\bar{w}) = iL^{0,\,1}_{\bar{w}}L^{1,\,0}_v(f)$. 

\end{Lem}

\noindent {\it Proof.} Both formulae being local, we can work in coordinates. Let $v=\sum_l v_l\,(\partial/\partial z_l)$ and $w=\sum_r w_r\,(\partial/\partial z_r)$, with the $v_l$'s and the $w_r$'s locally defined smooth functions. Using (\ref{eqn:iddbar-f_2-contraction_Lie_coord}), we get:  \begin{eqnarray}\label{eqn:iddbar-f_2-contraction_Lie_proof_1}(i\partial\bar\partial f)(v,\,\bar{w}) = i\sum\limits_{l,\,r}v_l\,\bar{w}_r\,L^{1,\,0}_{\frac{\partial}{\partial z_l}}L^{0,\,1}_{\frac{\partial}{\partial\bar{z}_r}}(f).\end{eqnarray}  

On the other hand, using (vii) of Lemma \ref{Lem:Lie-deriv-prop_0-1}, we get: \begin{eqnarray*}iL^{1,\,0}_vL^{0,\,1}_{\bar{w}}(f) = i\sum\limits_{l,\,r}v_l\,L^{1,\,0}_{\frac{\partial}{\partial z_l}}\bigg(\bar{w}_r\,L^{0,\,1}_{\frac{\partial}{\partial\bar{z}_r}}(f)\bigg) = i\sum\limits_{l,\,r}v_l\,\bar{w}_r\,L^{1,\,0}_{\frac{\partial}{\partial z_l}}L^{0,\,1}_{\frac{\partial}{\partial\bar{z}_r}}(f) + i\sum\limits_{l,\,r}v_l\,L^{1,\,0}_{\frac{\partial}{\partial z_l}}(\bar{w}_r)\,L^{0,\,1}_{\frac{\partial}{\partial\bar{z}_r}}(f).\end{eqnarray*} If $w$ is supposed holomorphic, $\frac{\partial\bar{w}_r}{\partial z_l} = 0$, hence $L^{1,\,0}_{\frac{\partial}{\partial z_l}}(\bar{w}_r)=0$, for all $l$ and $r$, so the last sum above vanishes. Comparing with (\ref{eqn:iddbar-f_2-contraction_Lie_proof_1}), this proves the contention of (i).

The equality in (ii) is proved in a similar way if $v$ is supposed holomorphic. \hfill $\Box$

\vspace{3ex}

\noindent {\it Proof of Theorem \ref{The:moment-map_def_final}.} $\bullet$ Let us fix an arbitrary element $(\xi_1,\dots , \xi_{n-2},\,\bar\eta_1,\dots , \bar\eta_{n-2})\in{\cal P}_\frg$ and holomorphic vector fields $v, w$ of type $(1,\,0)$ on ${\mathscr X}$. Thanks to Lemma \ref{Lem:eqn:iddbar-f_2-contraction_Lie}, on ${\mathscr X}$ we have:\begin{eqnarray}\label{eqn:moment-map_def_final_proof_1}\partial\bar\partial\bigg(\!\langle\mu(\cdot),\,(\xi_1,\dots , \xi_{n-2},\,\bar\eta_1,\dots , \bar\eta_{n-2})\rangle\!\bigg)(v,\,\bar{w}) = L^{1,\,0}_vL^{0,\,1}_{\bar{w}}\big(\langle\mu(\cdot),\,(\xi_1,\dots , \xi_{n-2},\,\bar\eta_1,\dots , \bar\eta_{n-2})\rangle\big).\end{eqnarray}

$\bullet$ Now, we fix a point $f\in{\mathscr X}$ (i.e. a holomorphic map $f:S\longrightarrow X$ such that $f^\star\omega_{n-1} = i\partial\bar\partial\Gamma$ on $S$ for some smooth real $(n-2,\,n-2)$-form $\Gamma$ on $S$). Consider a $2$-parameter family $(F_{s,\,t})_{(s,\,t)\in D(0,\,\varepsilon)\times D(0,\,\varepsilon)}$, where $D(0,\,\varepsilon)\subset\C$ is the open disc about $0$ of radius $\varepsilon$, of holomorphic maps $F_{s,\,t}:S\longrightarrow X$ such that $$F_{s,\,t}^\star\omega_{n-1} = i\partial\bar\partial\Gamma_{s,\,t}$$ on $S$ for all $(s,\,t)\in D(0,\,\varepsilon)\times D(0,\,\varepsilon)$ and some smooth real $(n-2,\,n-2)$-forms $\Gamma_{s,\,t}$ on $S$. Thus, $F_{s,\,t}\in{\mathscr X}$ for every $(s,\,t)$. We suppose, moreover, that: \begin{eqnarray*}(i)\,F_{0,\,0} = f;  \hspace{3ex} (ii)\,\frac{\partial F_{\cdot,\,0}}{\partial s}_{\bigg|s=0} = v(f)\in T^{1,\,0}_f{\mathscr X}; \hspace{3ex} (iii)\,\frac{\partial F_{0,\,\cdot}}{\partial t}_{\bigg|t=0} = w(f)\in T^{1,\,0}_f{\mathscr X}.\end{eqnarray*}

Equality (\ref{eqn:moment-map_def_final_proof_1}) and the definition of the Lie derivative in terms of the flow of the vector field w.r.t. which it is computed show that, on ${\mathscr X}$, we have the first equality below: \begin{eqnarray*}\nonumber\partial\bar\partial\bigg(\!\langle\mu(\cdot),\,(\xi_1,\dots , \xi_{n-2},\,\bar\eta_1,\dots , \bar\eta_{n-2})\rangle\!\bigg)(v,\,\bar{w}) = \frac{\partial^2}{\partial s\partial\bar{t}}_{\bigg|(s,\,t)=(0,\,0)}\bigg(\!\langle\mu(F_{s,\,\bar{t}}),\,(\xi_1,\dots , \xi_{n-2},\,\bar\eta_1,\dots , \bar\eta_{n-2})\rangle\!\bigg).\end{eqnarray*} Since \begin{eqnarray*}\langle\mu(F_{s,\,\bar{t}}),\,(\xi_1,\dots , \xi_{n-2},\,\bar\eta_1,\dots , \bar\eta_{n-2})\rangle = i\int\limits_S\Gamma_{s,\,\bar{t}}(\xi_1,\dots , \xi_{n-2},\,\bar\eta_1,\dots , \bar\eta_{n-2})\,dV,\end{eqnarray*} we conclude that 

\begin{eqnarray}\label{eqn:moment-map_def_final_proof_2}\partial\bar\partial\bigg(\!\!\langle\mu(\cdot),\,(\xi_1,\dots , \xi_{n-2},\,\bar\eta_1,\dots , \bar\eta_{n-2})\rangle\!\!\bigg)(v,\,\bar{w}) =\! \int\limits_S\!\!\bigg(i\frac{\partial^2\Gamma_{s,\,\bar{t}}}{\partial s\partial\bar{t}}_{\bigg|(s,\,t)=(0,\,0)}\bigg)(\xi_1,\dots , \xi_{n-2},\,\bar\eta_1,\dots , \bar\eta_{n-2})\,dV.\,\ \end{eqnarray}  

 $\bullet$ On the other hand, $v(f), w(f)\in T_f^{1,\,0}{\mathscr X} = H^0(S,\,f^\star T^{1,\,0}X)\subset H^0(S,\,T^{1,\,0}S)$. By using the definitions and the properties of $L^{1,\,0}$ and $L^{0,\,1}$ (see Definition \ref{Def:Lie-deriv_def} and Lemmas \ref{Lem:Lie-deriv-prop_1-0} and \ref{Lem:Lie-deriv-prop_0-1}), we get the following equalities on $S$: \begin{eqnarray*}\partial\bar\partial\bigg(\bar{w}(f)\lrcorner v(f)\lrcorner f^\star\omega_{n-1}\bigg) = \partial\bigg(L^{0,\,1}_{\bar{w}(f)}\bigg(v(f)\lrcorner f^\star\omega_{n-1}\bigg) - \bar{w}(f)\lrcorner\bar\partial\bigg(v(f)\lrcorner f^\star\omega_{n-1}\bigg)\bigg) \end{eqnarray*}
 \begin{eqnarray*}= \partial\bigg(\!v(f)\lrcorner L^{0,\,1}_{\bar{w}(f)}(f^\star\omega_{n-1}) + [\bar{w}(f),\,v(f)]\lrcorner f^\star\omega_{n-1}\!\bigg) - \partial\bigg(\!\bar{w}(f)\lrcorner(\bar\partial v(f))\lrcorner f^\star\omega_{n-1} - \bar{w}(f)\lrcorner v(f)\lrcorner\bar\partial(f^\star\omega_{n-1})\!\bigg).\end{eqnarray*} Now, $[\bar{w}(f),\,v(f)] = 0$ since $v(f)$ is a holomorphic $(1,\,0)$-vector field and $\bar{w}(f)$ is an anti-holomorphic $(0,\,1)$-vector field on $S$. Meanwhile, $\bar\partial v(f) = 0$ since $v(f)$ is holomorphic and $\bar\partial(f^\star\omega_{n-1}) = f^\star(\bar\partial\omega_{n-1}) = 0$ since $\omega$ is balanced.

 Thus, only the first of the four terms on the right of the above expression is non-vanishing, so we get the first of the following equalities on $S$: \begin{eqnarray*}\partial\bar\partial\bigg(\bar{w}(f)\lrcorner v(f)\lrcorner f^\star\omega_{n-1}\bigg) & \!\!=\! & \partial\bigg(\!v(f)\lrcorner L^{0,\,1}_{\bar{w}(f)}(f^\star\omega_{n-1})\!\bigg) = L^{1,\,0}_{v(f)} L^{0,\,1}_{\bar{w}(f)}(f^\star\omega_{n-1}) - v(f)\lrcorner\partial\bigg(L^{0,\,1}_{\bar{w}(f)}(f^\star\omega_{n-1})\bigg) \\
   & \!\!=\! & \frac{\partial^2}{\partial s\partial\bar{t}}_{\bigg|(s,\,t)=(0,\,0)}\bigg(F^\star_{s,\,\bar{t}}\,\omega_{n-1}\bigg) - v(f)\lrcorner L^{0,\,1}_{\bar{w}(f)}\bigg(\partial(f^\star\omega_{n-1})\bigg),\end{eqnarray*} where the last term is justified by property (iii)(b) of Lemma \ref{Lem:Lie-deriv-prop_0-1} and by the $(1,\,0)$-vector field $w(f)$ on $S$ being holomorphic.

 Since $\omega$ is balanced, we have $\partial(f^\star\omega_{n-1}) = f^\star(\partial\omega_{n-1}) = 0$, so the above equality reduces to the first of the following equalities on $S$: \begin{eqnarray*}\partial\bar\partial\bigg(\bar{w}(f)\lrcorner v(f)\lrcorner f^\star\omega_{n-1}\bigg) & \!= & \frac{\partial^2}{\partial s\partial\bar{t}}_{\bigg|(s,\,t)=(0,\,0)}\bigg(F^\star_{s,\,\bar{t}}\,\omega_{n-1}\bigg) = \frac{\partial^2}{\partial s\partial\bar{t}}_{\bigg|(s,\,t)=(0,\,0)}\bigg(i\partial\bar\partial\Gamma_{s,\,\bar{t}}\bigg) \\
 & \!= & \partial\bar\partial\bigg(i\frac{\partial^2\Gamma_{s,\,\bar{t}}}{\partial s\partial\bar{t}}_{\bigg|(s,\,t)=(0,\,0)}\bigg).\end{eqnarray*} Thus, we get the following equality on $S$: \begin{eqnarray}\label{eqn:moment-map_def_final_proof_3}i\,\frac{\partial^2\Gamma_{s,\,\bar{t}}}{\partial s\partial\bar{t}}_{\bigg|(s,\,t)=(0,\,0)} = \bar{w}(f)\lrcorner v(f)\lrcorner f^\star\omega_{n-1}.\end{eqnarray}

 $\bullet$ Putting together (\ref{eqn:moment-map_def_final_proof_2}) and (\ref{eqn:moment-map_def_final_proof_3}), we get: \begin{eqnarray*}\partial\bar\partial\bigg(\!\langle\mu(\cdot),\,(\xi_1,\dots , \xi_{n-2},\,\bar\eta_1,\dots , \bar\eta_{n-2})\rangle\!\bigg)(v,\,\bar{w}) =\!\int\limits_S\!\bigg(\!\bar{w}(f)\lrcorner v(f)\lrcorner f^\star\omega_{n-1}\!\bigg)(\xi_1,\dots , \xi_{n-2},\,\bar\eta_1,\dots , \bar\eta_{n-2})\,dV\end{eqnarray*} \begin{eqnarray*} & = & \int\limits_S(f^\star\omega_{n-1})\bigg(\xi_1,\dots , \xi_{n-2},\,\bar\eta_1,\dots , \bar\eta_{n-2},\,v(f),\,\bar{w}(f)\bigg)\,dV \\
     & = & \Omega\bigg(X(\xi_1),\dots , X(\xi_{n-2}),\,X(\bar\eta_1),\dots , X(\bar\eta_{n-2}),\,v,\,\bar{w}\bigg) \\
     & = & \bigg(X(\bar\eta_{n-2})\lrcorner\dots\lrcorner X(\bar\eta_1)\lrcorner X(\xi_{n-2})\lrcorner\dots\lrcorner X(\xi_1)\lrcorner\Omega\bigg)(v,\,\bar{w}),\end{eqnarray*} where the last-but-one equality follows from (\ref{eqn:Omega_def}). This proves (\ref{eqn:moment-map_def_final}) and completes the proof of Theorem \ref{The:moment-map_def_final}.  \hfill $\Box$

\vspace{3ex}

\noindent {\bf References.}

\vspace{1ex}

\noindent [ACPS23]\, D. Angella, S. Calamai, F. Pediconi, C. Spotti --- {\it A Moment Map for Twisted-Hamiltonian Vector Fields on Locally Conformally K\"ahler Manifolds} --- Transform. Groups  (2023).  DOI : 10.1007/s00031-023-09815-2

\vspace{1ex}

\noindent [AA87]\, L. Alessandrini, M. Andreatta --- {\it Closed Transverse $(p,\,p)$-forms on Compact Complex Manifolds} --- Compos. Math. {\bf 61} (1987), 181-200.

\vspace{1ex}

\noindent [Che87]\, P. Cherrier --- {\it \'Equations de Monge-Amp\`ere sur les vari\'et\'es hermitiennes compactes} --- Bull. Sc. Math. (2) {\bf 111} (1987), 343-385.

\vspace{1ex}

\noindent [Don99]\, S.K. Donaldson --- {\it Moment Maps and Diffeomorphisms} --- Asian J. Math. {\bf 3}, No. 1 (1999), 1-16.

\vspace{1ex}

\noindent [Fri89]\, R. Friedman --{\it On Threefolds with Trivial Canonical Bundle} --- in Complex Geometry and Lie Theory (Sundance, UT, 1989), Proc. Sympos. Pure Math., vol. 53, Amer. Math. Soc., Providence, RI, 1991, p. 103-134.

\vspace{1ex}

\noindent [GRT23]\, M. Garcia-Fernandez, R. Rubio, C. Tipler --- {\it Gauge Theory for String Algebroids} --- arXiv e-print DG 2004.11399v3, to appear in Journal of Differential Geometry. 

\vspace{1ex}

\noindent [GL09]\, B. Guan, Q. Li --- {\it Complex Monge-Amp\`ere Equations on Hermitian Manifolds} --- arXiv:0906.3548. 

\vspace{1ex}

\noindent [KP23]\, H. Kasuya, D. Popovici --- {\it Partially Hyperbolic Compact Complex Manifolds} ---  arXiv e-print DG 2304.01697v1.

\vspace{1ex}

\noindent [Kob70]\, S. Kobayashi --- {\it Hyperbolic Manifolds and Holomorphic Mappings} --- Marcel Dekker, New York (1970).

\vspace{1ex}

\noindent [MP22a]\, S. Marouani, D. Popovici --- {\it Balanced Hyperbolic and Divisorially Hyperbolic Compact Complex Manifolds} --- arXiv e-print CV 2107.08972v2, to appear in Mathematical Research Letters.

\vspace{1ex}

\noindent [MP22b]\, S. Marouani, D. Popovici --- {\it Some Properties of Balanced Hyperbolic Compact Complex Manifolds} --- Internat. J. Math., {\bf 33}, No. 3 (2022) 2250019, DOI : 10.1142/S0129167X22500197.

\vspace{1ex}

\noindent [Pop15]\, D. Popovici --- {\it Aeppli Cohomology Classes Associated with Gauduchon Metrics on Compact Complex Manifolds} --- Bull. Soc. Math. France {\bf 143} (3), (2015), p. 1-37.

\vspace{1ex}

\noindent [Sar78]\, K.S. Sarkaria --- {\it A Finiteness Theorem for Foliated Manifolds} --- J. Math. Soc. Japan {\bf 30}, No.4 (1978), 687-696.

\vspace{1ex}

\noindent [TW10]\, V. Tosatti, B. Weinkove --- {\it The Complex Monge-Amp\`ere Equation on Compact Hermitian Manifolds} --- J. Amer. Math. Soc. {\bf 23} (2010), no. 4, 1187-1195.

\vspace{1ex}

\noindent [Yac98]\, A. Yachou --- {\it Sur les vari\'et\'es semi-k\"ahl\'eriennes} --- PhD Thesis, University of Lille.

\vspace{1ex}

\noindent [Yau78]\, S.T. Yau --- {\it On the Ricci Curvature of a Complex K\"ahler Manifold and the Complex Monge-Amp\`ere Equation I} --- Comm. Pure Appl. Math. {\bf 31} (1978) 339-411.

\vspace{3ex}

\noindent Institut de Math\'ematiques de Toulouse, Universit\'e Paul Sabatier,

\noindent 118 route de Narbonne, 31062 Toulouse, France

\noindent Email: popovici@math.univ-toulouse.fr

\vspace{2ex}

\noindent and

\vspace{2ex}

\noindent Departamento de Matem\'aticas\,-\,I.U.M.A., Universidad de Zaragoza,

\noindent Campus Plaza San Francisco, 50009 Zaragoza, Spain

\noindent Email: ugarte@unizar.es

\end{document}